\documentclass{artjlt}

\usepackage{amsmath,amssymb,amsfonts,amsxtra}
\usepackage{paralist}
\usepackage{dsfont,mathrsfs}
\renewcommand{\mathcal}{\mathscr}
\usepackage[all]{xypic}
\xyoption{dvips}

\newcommand{\A}{\mathrm{A}}
\newcommand{\beq}{\begin{equation}}
\newcommand{\C}{\mathscr{C}}
\newcommand{\cc}{\mathfrak{c}}
\newcommand{\CC}{\mathds{C}}
\newcommand{\dl}{\mathfrak{l}}
\newcommand{\eeq}{\end{equation}}
\newcommand{\eps}{\epsilon}
\newcommand{\gl}{\mathfrak{gl}}
\newcommand{\sll}{\mathfrak{sl}}
\newcommand{\dd}{\mathfrak{d}}
\newcommand{\fg}{\mathfrak{g}}
\newcommand{\hh}{\mathfrak{h}}
\newcommand{\rr}{\mathfrak{r}}
\newcommand{\fs}{\mathfrak{s}}
\newcommand{\ideal}[1]{{\left\langle#1\right\rangle}}
\newcommand{\il}{\overline{\mathfrak{l}}}

\newcommand{\mm}{\mathfrak{m}}
\newcommand{\nn}{\mathfrak{n}}
\newcommand{\NN}{\mathds{N}}
\newcommand{\sO}{\mathscr{O}}
\newcommand{\ol}{\overline}
\newcommand{\ul}{\underline}
\newcommand{\onto}{\twoheadrightarrow}
\newcommand{\p}{\partial}
\newcommand{\sP}{\mathscr{P}}
\newcommand{\QQ}{\mathds{Q}}

\newcommand{\wh}{\widehat}
\newcommand{\xymat}{\SelectTips{cm}{}\xymatrix}

\DeclareMathOperator{\Ad}{Ad}
\DeclareMathOperator{\Aut}{Aut}
\DeclareMathOperator{\Der}{Der}
\DeclareMathOperator{\End}{End}
\DeclareMathOperator{\GL}{GL}
\DeclareMathOperator{\ord}{ord}
\DeclareMathOperator{\rk}{rk}
\DeclareMathOperator{\Sing}{Sing}
\DeclareMathOperator{\wt}{wt}

\numberwithin{equation}{section}

\title{Initial logarithmic Lie algebras of\\ hypersurface singularities}
\author{Michel Granger and Mathias Schulze\thanks{MS was supported by the College of Arts \& Sciences at Oklahoma State University through a FY 2009 Dean's Incentive Grant.}}
\lastname{Granger and Schulze}



\address{Michel Granger\\
D\'epartement de Mathematiques\\
Universit\'e d'Angers\\
2 Bd Lavoisier\\
49045 Angers\\
France\\
granger@univ-angers.fr}

\address{Mathias Schulze\\
Oklahoma State University\\
401 Mathematical Sciences\\
Stillwater, OK 74078\\
United States\\
mschulze@math.okstate.edu}

\begin{document}

\maketitle

\begin{abstract}
We introduce a Lie algebra of initial terms of logarithmic vector fields along a hypersurface singularity.
Extending the formal structure theorem in \cite[Thm.~5.4]{GS06}, we show that the completely reducible part of its linear projection lifts formally to a linear Lie algebra of logarithmic vector fields.
For quasihomogeneous singularities, we prove convergence of this linearization.
We relate our construction to the work of Hauser and M\"uller \cite{Mue86,HM89} on Levi subgroups of automorphism groups of singularities, which proves convergence even for algebraic singularities.

Based on the initial Lie algebra, we introduce a notion of reductive hypersurface singularity and show that any reductive free divisor is linear.

As an application, we describe a lower bound for the dimension of hypersurface singularities in terms of the semisimple part of their initial Lie algebra.
\end{abstract}

\section{Initial logarithmic Lie algebras}\label{21}

Let $X=(\CC^n,0)$ the space germ of $\CC^n$, $n\ge1$, at the origin with local ring $\sO=\sO_X$.
A (reduced) function $f\in\mm\subset\sO$ where $\mm$ is the maximal ideal of $\sO$ defines a hypersurface singularity 
\[
0\in D=\{f=0\}\subset X
\]
Let $\Der_\CC(\sO)$ be the $\sO$-module of germs of vector fields at $0$ and 
\[
\Der(-\log D)=\{\theta\in\Der_\CC(\sO)\mid\theta(f)\in\sO\cdot f\}
\]
the $\sO$-module of logarithmic vector fields along $D$ introduced by K.~Saito \cite{Sai80}.
We shall always assume that
\beq\label{1}
\dl_D:=\Der(-\log D)\subset\mm\cdot\Der_\CC(\sO)=:\Delta
\eeq
which means by Rossi's theorem \cite[Cor.~3.4]{Ros63} that $D$ is not isomorphic to a product $D'\times\CC$ with a smooth factor for some smaller dimensional $D'\subset(\CC^{n-1},0)$.
Note that 
\beq\label{2}
[\mm^k\cdot\Delta,\mm^l\cdot\Delta]\subset\mm^{k+l}\cdot\Delta
\eeq
and hence
\[
\Delta^k:=\mm^k\cdot\Delta,\quad\Delta_k:=\Delta/\Delta^{k+1},\quad k\ge0,
\]
are Lie algebras.
Note that $\Delta_0=\End_\CC(\mm/\mm^2)$.
There are maps of Lie algebras
\beq\label{3}
\xymat{\Delta\ar@{->>}[r]^-{\pi_k}&\Delta_k\ar@{->>}[r]^-{\pi^k_l}&\Delta_l},\quad k>l\ge0.
\eeq

\begin{Definition}\label{25}
We call 
\[
\il_D:=\dl_D/\mm\cdot\dl_D
\]
the \emph{initial logarithmic Lie algebra} of $D$.
\end{Definition}

For $\delta,\theta\in\dl_D$ and $p\in\mm$, $\delta(p)\in\mm$ by assumption \eqref{1} and hence
\[
[\delta,p\cdot\theta]=\delta(p)\cdot\theta+p\cdot[\delta,\theta]\in\mm\cdot\dl_D.
\]
This shows that $\il_D$ is a Lie algebra and $\pi_0$ in \eqref{3} induces
\beq\label{4}
\xymat{\dl_D\ar@{->>}[r]^-{\ul\lambda}&\il_D\ar@{->>}[r]^-{\ol\lambda_0}&\dl_{D,0}:=\pi_0(\dl_D)\subset\Delta_0}
\eeq
with
\beq\label{40}
\ker\ol\lambda_0=(\Delta^1\cap\dl_D)/\mm\cdot\dl_D.
\eeq

\begin{Example}\label{5}
For a normal crossing divisor $D=\{x_1\cdots x_n=0\}$,
\[
\il_D=\sum_{i=1}^n\CC\cdot x_i\p_i\cong\CC^n.
\]
\end{Example}

\begin{Remark}
In \cite[Prop.~6.2]{GS06} we showed that $\dl_D/(\dl_D\cap\mm^k\cdot\Delta)$ is solvable for all $k\ge0$ if $D$ is a free divisor in dimension $n\le3$.
By Krull's intersection theorem, the denominator is contained in $\mm\cdot\dl_D$ for large $k$.
This implies that $\il_D$ is solvable in this case.
\end{Remark}

In Section~\ref{22}, we shall relate the size of the singular locus $\Sing D$ of $D$ to the size of a Levi factor of $\il_D$.
The following proposition proved in Section~\ref{22} and the following example taken from \cite[\S 6.4]{GMN06} serve as a motivation for this plan.

\begin{Proposition}\label{6}
If $D$ is an isolated singularity of order at least $3$ then $\il_D$ is solvable.
\end{Proposition}

\begin{Remark}
Also for $n=2$ the conclusion of Proposition~\ref{6} holds true as $D$ is a free divisor and hence $\dim\il_D=n<3=\dim\sll_2(\CC)$.
However the statement of Proposition~\ref{6} is in general wrong for $\ord D=2$ and $n\ge3$.
Consider for example $D$ defined by $f=x_1^2+\cdots+x_n^2$ where $\il_D$ is generated by $\eps=x_1\p_1+\cdots+x_n\p_n$ and $x_i\p_j-x_j\p_i$ with $1\le i<j\le n$, see \eqref{18} and \eqref{19}.
For $n\ge3$, the latter elements span the nonsolvable Lie algebra $o_n(\CC)$ of all skew symmetric complex square matrices.
\end{Remark}

\begin{Example}\label{7}
The divisor $D\subset\CC^4$ defined by
\[
f=y^2z^2-4xz^3-4y^3w+18xyzw-27w^2x^2=0
\]
has a basis of logarithmic vector fields 
\beq\label{8}
\begin{pmatrix}
x&y&z&w\\
-3x&-y&z&3w\\
y&2z&3w&0\\
0&3x&2y&z
\end{pmatrix}
\cdot
\begin{pmatrix}
\p_x\\\p_y\\\p_z\\\p_w\\
\end{pmatrix}.
\eeq
Hence $\il_D\cong\gl_2(\CC)=\sll_2(\CC)\oplus\CC$ and $\mm/\mm^2$ is an irreducible $\sll_2(\CC)$-representation.
This is an example of a free divisor which means that $\dl_D$ is a free $\sO$-module.
By Saito's criterion \cite[Lem.~1.8.ii]{Sai80}, this follows from the determinant of the left matrix in \eqref{8} being a unit multiple of $f$.
The singular locus of a free divisors is equidimensional of the maximal possible dimension.
In our case,
\beq\label{9}
\Sing S=\{z^2-3yw=yz-9xw=y^2-3xz=0\},\quad\dim\Sing D=2.
\eeq
\end{Example}

In Theorem~\ref{13}, the preceding observations will be generalized to a lower bound for $\Sing D$ in terms of the weight diagram of a Levi factor of $\il_D$ acting on $\mm/\mm^2$ via $\ol\lambda_0$.
Denoting completion at $0$ and $\mm$-adic completion by $\wh\ $, we have
\[
\Der(-\log\wh D)=\wh\Der(-\log D),\quad\il_D=\il_{\wh D},
\] 
by exactness of completion and by Definition~\ref{25}.
By faithful flatness of completion, we also have
\[
\dim\Sing D=\dim\Sing\wh D.
\]
Therefore our problem is purely formal.
In Section~\ref{24}, we shall develop the main technical tool for our investigations which is an extension of the formal structure theorem for $\Der(-\log D)$ from \cite[Thm.~5.4]{GS06}.
Theorem~\ref{17} states that any completely reducible Lie subalgebra of $\dl_{D,0}$ can be lifted to $\wh\dl_D$ and linearized in formal suitable coordinates.
In loc.\ cit.\ this was only proved for Abelian Lie subalgebras of semisimple vector fields.
The proof of this extension is essentially based on the idea of \cite{Her68}.

In Section~\ref{31}, we introduce the notion of a reductive hypersurface singularity by essentially requiring that $\il_D\cong\dl_{D,0}$ is reductive.
This is motivated by the recent interest in reductive linear free divisors \cite{BM06,GMN06,GS08,GMS08}, Example~\ref{7} being a representative of this class.
We shall show in Theorem~\ref{32} that a reductive free divisor is automatically linear.
More generally, we show in Theorem~\ref{39} that our formal linearization of the completely reducible part of $\dl_{D,0}$ in Theorem~\ref{17} is convergent if $D$ is a quasihomogeneous singularity.
The proof of these convergence results is based on Artin's Approximation Theorem \cite{Art68}.
 
Finally, we show that our linearized completely reducible Lie algebra is the Lie algebra of a Levi subgroup of the embedded automorphism group of the singularity.
Using results of Hauser and M\"uller \cite{Mue86,HM89}, this yields convergence in Theorem~\ref{17} even for algebraic singularities, as formulated in Theorem~\ref{46}.

\section{Linearization of completely reducible initial Lie subalgebras}\label{24}

For a Lie algebra $\fg$ with solvable radical $\rr\subsetneq\fg$ the natural projection $\fg\onto\fg/\rr$ has a section by Levi's theorem \cite[Ch.~III, \S9]{Jac79}. 
Its semisimple image $\fs$ in $\fg$ is called a Levi factor of $\fg$ and $\fg=\fs\oplus\rr$ is called a Levi decomposition.

We call a Lie subalgebra $\fg$ of $\wh\Delta$ formally solvable/nilpotent if $\pi_k(\fg)$ is solvable/nilpotent for all $k\in\NN$. 
Then any Lie subalgebra $\fg$ of $\wh\Delta$ contains a formal solvable radical $\rr$ such that $\fg/\rr$ is a finite semisimple Lie algebra.
Indeed, as $\Delta^1=\ker\pi_0$ is formally nilpotent, $\rr$ is just the preimage of the solvable radical of $\pi_0(\fg)\subset\Delta_0$ under the restriction of $\pi_0$ to $\fg$, see the proof of Theorem~\ref{17}.\eqref{17a}.
If $\fs$ is the image of a section of the natural projection $\fg\onto\fg/\rr$, we call $\fg=\fs\oplus\rr$ a formal Levi decomposition of $\fg$.
In the same way, we define $\hh\subset\fg$ to be a formal Cartan subalgebra of $\fg$ if $\pi_k(\hh)\subset\pi_k(\fg)$ is a Cartan subalgebra of $\pi_k(\fg)$ in the sense of \cite[Ch.~III, \S1, Def.~1]{Ser01} for all $k\in\NN$.
This means that $\hh$ is formally nilpotent and equals its normalizer in $\fg$.

A basis of $\mm$ corresponds to a section $\phi$ of $\pi_0$ in \eqref{3}.
For $\delta\in\Delta$, we call $\phi$ of the semisimple part of the endomorphism $\pi_0(\delta)\in\Delta_0$ the semisimple part of $\delta$ with respect to the basis corresponding to $\phi$.

\begin{Theorem}\label{17}
Let $\ol\dl=\il_D$ and $\dl=\wh\dl_D$ with solvable radicals $\ol\rr$ and $\rr$.
Then the following holds:
\begin{enumerate}[(a)]

\item\label{17a} Any Levi decomposition $\ol\dl=\ol\fs\oplus\ol\rr$ lifts to a formal Levi decomposition $\dl=\fs\oplus\rr$ and $\fs$ is linearizable.

\item\label{17b} The solvable radical $\rr_0=\pi_0(\rr)$ of the linear part $\dl_0=\pi_0(\dl)$ of $\dl_0$ can be decomposed as $\rr_0=\dd_0\oplus\nn_0$ where $\dd_0$ is a Lie subalgebra with semisimple elements that commutes with the Levi factor $\fs_0=\pi_0(\fs)$, and $\nn_0$ is a nilpotent ideal.

\item\label{17c} With respect to a basis $\wh\phi$ of $\wh\mm$, $\fs_0$ lifts to a $\wh\phi$-linear formal Levi factor of $\dl$, and $\dd_0$ lifts to a $\wh\phi$-linear Abelian Lie subalgebra $\dd\subset\rr$ with semisimple elements that commutes with $\fs$.
In particular, $\rr=\dd\oplus\nn$ where $\nn=\pi_0^{-1}(\nn_0)\cap\dl$.

\item\label{17d} Let $\hh$ be a Cartan subalgebra of $\fs$ and let $\cc$ be the centralizer of $\hh\oplus\dd$ in $\dl$.
Then $\cc$ is a formal Cartan subalgebra of $\dl$ and $\hh\oplus\dd$ is the Lie subalgebra of semisimple parts of elements in $\cc$ (with respect to $\wh\phi$ in \eqref{17c});

\item\label{17e} The representation of $\fs\oplus\dd$ in $\End_\CC(\mm/\mm^2)$ is defined over $\QQ$ (with respect to $\wh\phi$ in \eqref{17c}).

\item\label{17f} $\wh D$ is defined by a $(\hh\oplus\dd)$-homogeneous $f\in\wh\mm$. 
Furthermore $\fs$ and $\rr$ are both $\hh$- and $\dd$-graded modules with trivial $\dd$-graduation on $\fs$.

\end{enumerate}
\end{Theorem}

\begin{Proof}
\begin{asparaenum}[(a)]
\item Consider the Lie algebras
\[
\dl_k:=\pi_k(\dl)\cong\dl/\dl\cap\Delta^{k+1},\quad k\ge0,
\]
with solvable radical $\rr_k\subset\dl_k$ and projections
\[
\xymat{\dl\ar@{->>}[r]^-{\lambda_k}&\dl_k\ar@{->>}[r]^-{\lambda^k_l}&\dl_l,\quad k>l\ge0}
\]
induced by \eqref{3}.
By Krull's intersection theorem there is an $l\ge0$ such that
\beq\label{41}
\Delta^{l+1}\cap\dl\subset\mm\cdot\dl
\eeq
and hence a projection
\[
\xymat{\dl_l\ar@{->>}[r]^-{\ul\lambda^l}&\il}.
\]
As $\ker\ul\lambda^l=(\mm\cdot\dl+\Delta^{l+1})/\Delta^{l+1}\subset\Delta^1/\Delta^{l+1}$ and $\ker\lambda^k_l=\ker\pi^k_l\cap\dl_k\subset\Delta^1/\Delta^{k+1}$ are nilpotent by \eqref{2} and hence solvable, we have that $(\ul\lambda^l)^{-1}(\ol\rr)=\rr_l$ and $(\lambda^k_l)^{-1}(\rr_l)=\rr_k$ by \cite[Lem.~I.7]{Jac79}.
We can then construct lifted Levi decompositions
\beq\label{27}
\dl_k=\fs_k\oplus\rr_k,\quad\dl_l=\fs_l\oplus\rr_l,\quad\xymat{\fs_k\ar[r]^-{\lambda^k_l|_{\fs_k}}_-\cong&\fs_l\ar[r]^-{\ul\lambda^l|_{\fs_l}}_-\cong&\ol\fs}.
\eeq
We shall explain how to lift $\dl_l=\fs_l\oplus\rr_l$ to $\dl_k=\fs_k\oplus\rr_k$ along $\lambda^k_l$.
The former is obtained analogously by lifting $\ol\dl=\ol\fs\oplus\ol\rr$ along $\ul\lambda^l$.
The map $\lambda^k_l$ induces a projection
\[
\xymat@C=50pt{(\lambda^k_l)^{-1}(\fs_l)\ar@{->>}[r]^-{\lambda^k_l|_{(\lambda^k_l)^{-1}(\fs_l)}}&\fs_l}
\]
whose nilpotent kernel $\ker\lambda^k_l$ must be the solvable radical of $(\lambda^k_l)^{-1}(\fs_l)$.
By Levi's Theorem \cite[Ch.~III, \S9]{Jac79}, there is a Levi decomposition
\[
(\lambda^k_l)^{-1}(\fs_l)=\fs_k\oplus\ker\lambda^k_l
\]
such that $\lambda^k_l$ induces the isomorphism in \eqref{27}.
Taking the projective limits $\fs=\varprojlim_k\fs_k $ and setting $\rr:=\ul\lambda^{-1}(\ol\rr)$, we obtain a formal lifting
\[
\dl=\fs\oplus\rr,\quad\xymat{\fs\ar[r]^-{\ul\lambda|_\fs}_-\cong&\ol\fs}
\]
of the Levi decomposition $\ol\dl=\ol\fs\oplus\ol\rr$.
We set  $\fs_k:=\lambda_k(\fs)$ and $\rr_k:=\lambda_k(\rr)$ for all $k\ge0$.
Because $\ol\rr=\ol\lambda_0^{-1}(\rr_0)$ and $\rr_k=(\lambda^k_0)^{-1}(\rr_0)$, there are Levi decompositions
\[
\dl_k=\fs_k\oplus\rr_k,\quad\xymat{\ol\fs\ar[r]^-{\ol\lambda_0|_{\ol\fs}}_-\cong&\fs_0&\ar[l]_-{\lambda^k_0|_{\fs_k}}^-\cong\fs_k}.
\]
By \cite[Thm.~III.1]{Her68}, $\fs$ can be linearized.

\item By the Theorem of Malcev--Harish-Chandra \cite[Ch.~III, \S9]{Jac79} all Levi factors are conjugate by an inner automorphism, so we may ignore the requirement $\pi(\fs)=\fs_0$.
By \cite[Lem.~1]{Sch06}, $\dl_0$ is the Lie algebra of a linear algebraic group $A_0$.
Now the theorem of Mostow~\cite{Mos56} on the existence of a Levi decomposition of this group yields the claim.
However we present a more elementary argument:

As $\rr_0$ is the Lie algebra of a linear algebraic group, namely that of the solvable radical of $A_0$, \cite[\S15.2]{Hum75} shows that $\rr_0$ is closed under the operation of taking the semisimple or nilpotent part of an element.
Moreover, $\rr_0$ can be triangularized by Lie's Theorem \cite[Ch.~II,\S6]{Jac79} and then the strictly triangular or nilpotent elements of $\rr_0$ form an ideal $\nn_0\subset\rr_0$.
Let $\dd_0$ be a maximal dimensional Lie subalgebra of $\rr_0$ consisting of diagonal matrices in some coordinate system that triangularizes $\rr_0$.
With respect to $\dd_0$, $\dl_0$ decomposes as a direct sum of weight spaces
\beq\label{43}
\dl_0^\chi=\{\delta\in\dl_0\mid\forall\sigma\in\dd_0\colon[\sigma,\delta]=\chi(\sigma)\cdot\delta\},\quad\chi\in\dd_0^*,
\eeq 
corresponding to a block decomposition of matrices.
Also $\rr_0$ decomposes as a direct sum of weight spaces $\rr_0^\chi=\dd_0^\chi\cap\rr_0$, $\chi\in\dd_0^*$.
The Lie subalgebras $\dl_0^0\subset\dl_0$ and $\rr_0^0\subset\rr_0$ consists exactly of the block-diagonal elements.
Denote $\rr_0^{\ne0}=\bigoplus_{0\ne\chi\in\dd_0^*}\rr_0^\chi$ the complement of $\rr_0$ consisting of strictly block-triangular elements.
Modulo $\nn_0$, any element of $\delta\in\rr_0$ lies in $\rr_0^0$ and is semisimple.
So it can be diagonalized keeping $\rr_0$ triangular and $\dd_0$ diagonal.
By the maximality assumption, this means that $\delta\in\dd_0$.
We conclude that $\rr_0=\dd_0\oplus\nn_0$. 

Now let $\fs'_0$ be any Levi factor of $\dl_0$.
As $\rr_0$ is an ideal and $\dd_0$ acts by characters on the matrix blocks, $[\fs'_0,\dd_0]\subset\rr_0^{\ne0}\cap\nn_0$ and hence $\fs'_0\subset\dl_0^0\oplus(\rr_0^{\ne0}\cap\nn_0)$ consists of block-diagonal matrices.
Since $\fs'_0\cap\rr_0=0$, the projection onto the $\dl_0^0$ is a monomorphism of Lie algebras.
Its image $\fs_0\cong\fs'_0$ is a Levi factor of $\dl_0$ that commutes with $\dd_0$.

\item First, we construct a linear lift $\dd\subset\rr$ of $\dd_0$ with respect to some basis of $\wh\mm$.
To this end, pick $\delta\in\rr$ such that $\lambda_0(\delta)\in\dd_0$ and consider its Poincar\'e--Dulac decomposition \cite[Ch.~3.~\S3.2]{AA88}
\beq\label{28}
\delta=\sigma+\nu,\quad[\sigma,\nu]=0,
\eeq
into commuting semisimple and nilpotent part, with respect to some basis of $\wh\mm$.
The proof of \cite[Thm.~5.4]{GS06} shows that the decomposition \eqref{28} takes place in $\dl$, in particular, $\nu\in\dl$.
On the other hand, $\lambda_0(\delta)=\lambda_0(\sigma)+\lambda_0(\nu)\in\dd_0$ is semisimple and hence $\nu\in\Delta^1$.
But by nilpotency of $\Delta^1$, we have $\nu\in\Delta^1\cap\dl\subset\rr$ and hence $\sigma\in\rr$.
As $\lambda_0(\delta)=\lambda_0(\sigma)$, we may therefore assume that $\delta=\sigma$ is linear.
By \cite[Thm.~5.3]{GS06}, the Poincar\'e--Dulac decomposition \eqref{28} can be iterated to simultaneously lift all elements in $\dd_0$ to a linear Abelian Lie subalgebra $\dd\subset\rr$ of semisimple elements.
Let $\dl^\chi$ denote the weight space with respect to $\chi\in\dd^*$ defined as in \eqref{43}.
The construction from \eqref{17a} applied to the centralizer $\dl^0\subset\dl$ of $\dd$ serves to lift $\fs_0\subset\dl_0^0$ to $\dl^0$.

Now, a slight modification of construction in the proof of \cite[Thm.~III.1]{Her68} serves to linearize $\fs$ without destroying the linearity of $\dd$.
To see this, replace in loc.\ cit.\ the space $V^j$ of homogeneous elements of degree $j$ in $\wh\Delta$ by its intersection $V^j_0$ with the centralizer of $\dd$ in $\wh\Delta$.
As $\fs$ and $\dd$ commute, $\fs$ is contained in this centralizer and each $V^j_0$ is an $\fs$-module via $\lambda_0$.
Then the infinitesimal coordinate changes $W^j\in\wh\Delta$ linearizing $\fs$ constructed in loc.\ cit.\ are in the centralizer of $\dd$ in $\wh\Delta$.
By \cite[\S II]{Her68}, the resulting infinitesimal coordinate change $W$ yields a $\dd$-homogeneous coordinate change $\exp(\Ad W)$ that leaves $\dd$ invariant by \cite[Lem.~2.7]{GS06}.
This gives the desired basis $\wh\phi$ of $\wh\mm$.

\item By \cite[Ch.~III, \S5, Thm.~3.(b)]{Ser01}, $\hh$ equals its centralizer in $\fs$. 
So $\cc=\hh\oplus\dd'$, where $\dd'$ is the centralizer of $\hh\oplus\dd$ in $\rr$.
For $\chi\in(\hh\oplus\dd)^*$, consider the $\chi$-weight spaces $\rr^\chi$ defined as in \eqref{43}.
As $\nn=\lambda^{-1}(\nn_0)\subset\rr$ is an ideal, it is $(\hh\oplus\dd)$-homogeneous and $\nn^\chi=\nn\cap\rr^\chi$.
Then $\dd'=\rr^0$ is a Lie subalgebra of $\rr$ and 
\beq\label{35}
\dd'=\dd\oplus\nn^0,\quad\nn^0=\nn\cap\dd',
\eeq
is a commuting sum of the Abelian Lie subalgebra $\dd\subset\dd'$ and the nilpotent ideal $\nn^0\subset\dd'$.
Thus, $\dd'$ is formally nilpotent and its semisimple parts lie in $\dd$.

It remains to show that $\cc$ equals its normalizer in $\hh\oplus\rr$ which reduces to prove that if $\eta\in\rr$ normalizes $\hh\oplus\dd'$ then $\eta\in\dd'$. 
Expanding it as
\[
\eta=\sum_{\chi\in(\hh\oplus\dd)^*}\eta_\chi,\quad\eta_\chi\in\rr^\chi,
\]
we obtain that
\[
\forall\delta\in\dd\subset\dd'\colon[\delta,\eta]=\sum_{\chi\in\dd^*}\chi(\delta)\cdot\eta_\chi\in\dd'=\rr^0.
\]
Hence, $\eta_\chi=0$ for $\chi\ne0$ and therefore $\eta=\eta_0\in\rr^0=\dd'$.

\item By complete reducibility \cite[Ch.~III, \S7, Thm.~8]{Jac79}, the $\fs$-module $V=\mm/\mm^2=\oplus_{j=1}^mV_j$ decomposes into irreducible $\fs$-modules $V_j$.
Each of these modules is defined over $\QQ$ by semisimplicity of $\fs$ \cite[Lem.~2.6]{GS08}.
As $\fs$ and $\dd$ commute and $\dd$ is Abelian, $\dd$ acts on each $V_j$ through a character $\chi_j\in\dd^*$, that is, $\delta(v_j)=\chi_j(\delta)\cdot v_j$ for $\delta\in\dd$ and $v_j\in V_j$.
By the construction in \cite[Lem.~1.4]{Sai71}, there is a basis of $\dd$ on which these characters take values in $\QQ$.

\item A $(\hh\oplus\dd)$-homogeneous $f\in\wh\mm$ defining $\wh D$ is constructed as in the proof of \cite[Thm.~5.4]{GS06}.
The second statement is clear since $\hh\subset\fs$, $\dd$ commutes with $\fs$, and $\rr$ is an ideal in $\dl$.
\end{asparaenum}
\end{Proof}

\begin{Remark}\label{30}
\begin{asparaenum}
\item\label{30a} Theorem~\ref{17}.\eqref{17d} shows that the space $\hh\oplus\dd$ satisfies the maximality condition in \cite[Thm.~5.4.(2)]{GS06}. 
So $s_D:=\dim_\CC(\hh\oplus\dd)$ is the maximal multihomogeneity of $D$ defined in \cite{Sch06}.
\item\label{30b} By \cite[Ch.~III, \S4, Cor.~1]{Ser01}, the rank $\rk\dl_0$ of the Lie algebra $\dl_0$ is the dimension 
of any of its Cartan subalgebras. 
According to \eqref{35} it is equal to  
\begin{align*}
\rk\dl_0=\dim_\CC\lambda_0(\dd')
&=\dim_\CC(\hh\oplus\dd)+\dim_\CC(\nn\cap\lambda_0(\dd'))\\
&=\dim_\CC(\hh\oplus\dd)+\dim_\CC(\nn\cap\lambda_0(\cc)).
\end{align*}
Here $n_D:=\dim_\CC(\nn\cap\lambda_0(\cc))$ is the dimension of the Lie algebra of nilpotent elements in any Cartan subalgebra of $\dl_0$.
This shows that $s_D$ is intrinsically defined as
\[
s_D=\rk\dl_0-n_D.
\]
\end{asparaenum}
\end{Remark}

\section{Reductive free divisors and algebraic singularities}\label{31}

\begin{Definition}\label{42}
\begin{enumerate}[(a)]
\item\label{42a} We call $D$ and $\wh D$ \emph{reductive} if $\ol\lambda_0|_{\ol\rr}\colon\ol\rr\to\dd_0$ in Theorem~\ref{17} is an isomorphism.
\item\label{42b} We call $D$ \emph{linear} if $\dl_D$ admits linear generators in some (convergent) coordinate system.
Similarly, we call $\wh D$ \emph{linear} if $\dl_{\wh D}$ admits linear generators in some (formal) coordinate system.
\end{enumerate}
\end{Definition}

With this terminology we draw the following conclusion from Theorem~\ref{17}.

\begin{Corollary}\label{37}
Any formal reductive hypersurface singularity $\wh D$ is linear.
\end{Corollary}

\begin{Proof}
By Nakayama's Lemma, Definition~\ref{42}.\eqref{42a} implies that the linearization $\fs\oplus\dd$ in Theorem~\ref{17} generates $\dl$.
\end{Proof}

Following K.~Saito, we call $D$ free if $\dl_D$ is a free $\sO$-module.
Then Saito's Criterion \cite[Thm.~1.8.ii]{Sai80} states that, for $\delta_1,\dots,\delta_n\in\dl_D$, $D$ is free with $\delta_1,\dots,\delta_n$ being a basis of $\dl_D$ if and only if
\beq\label{33}
f=\det(\delta_i(x_j))
\eeq
is a reduced defining equation of $D$.
In case of a formal hypersurface singularity $\wh D$, the latter property implies the former one \cite[Prop.~4.2]{GS06}.

\begin{Proposition}\label{34}
Let $D$ be a free divisor and let $\wh\phi\in\Aut(\wh X)$ with inverse $\wh\psi$ such that $\wh\phi_*\wh\dl_D=\dl_{\wh\phi(\wh D)}$ (which is $\wh\sO$-free) has a (convergent) $\wh\sO$-basis $\delta'_1,\dots,\delta'_n\in\wh\phi_*\Delta$.
Then $f'$ defined as in \eqref{33} using the basis $\delta'_1,\dots,\delta'_n$ and $x'_j=\wh\psi^*x_j$ defines a (free) divisor $D'$ isomorphic to $D$ and $\delta'_1,\dots,\delta'_n$ is a basis of $\dl_{D'}$.
\end{Proposition}

\begin{Proof}
Let $\delta_1,\dots,\delta_n$ be an $\sO$-basis of $\dl_D$.
Then $f$ in \eqref{33} is a reduced defining equation of $D$ and 
\[
(\delta_1,\dots,\delta_n)\cdot\wh U=(\wh\psi_*\delta'_1,\dots,\wh\psi_*\delta'_n)
\]
for some invertible matrix $\wh U\in\wh\sO^{n\times n}$.
Note that 
\[
\wh\psi_*\p_{x'_k}=\sum_{j=1}^n\frac{\p\wh\psi_j}{\p x'_k}\circ\wh\phi\cdot\p_{x_j}
\]
and hence
\[
(\wh\psi_*\p_{x'_k})(x_j)=\wh\phi^*\frac{\p\wh\psi_j}{\p x'_k}
\]
which implies that
\begin{align*}
\wh\psi_*(\delta'_i)(x_j)&=\sum_{k=1}^n\wh\psi_*(\delta'_i(x'_k)\p_{x'_k})(x_j)
=\sum_{k=1}^n\wh\psi_*(\delta'_i(x'_k))\cdot(\wh\psi_*\p_{x'_k})(x_j)\\
&=\sum_{k=1}^n\wh\phi^*(\delta'_i(x'_k))\cdot\wh\phi^*\frac{\p\wh\psi_j}{\p x'_k}.
\end{align*}
In matrix notation, these equalities read 
\[
((\wh\psi_*\delta'_i)(x_j))=(\wh\phi^*(\delta'_i(x'_j)))\circ\wh\phi^*(D\wh\psi)^t
\]
Setting $\wh u=\det(\wh U)/\det(\wh\phi^*(D\wh\psi))\in\wh\sO^*$, we conclude that
\beq\label{36}
\wh\phi^*f'=\det(\wh\phi^*(\delta'_i(x'_j)))=\wh u\cdot \det(\delta_i(x_j))=\wh u\cdot f.
\eeq
This shows that $\wh\phi(\wh D)=\wh D'$ and hence $\delta'_1\,\dots,\delta'_n\in\dl_{\wh D'}=\wh\dl_{D'}$.
But then convergence of $f'$ and $\delta'_1\,\dots,\delta'_n$ immediately implies that $\delta'_1\,\dots,\delta'_n\in\dl_{D'}$.
Due to faithful flatness of completion or Saito's Criterion, $\delta'_1,\dots,\delta'_n$ even form a basis of $\dl_{D'}$.

Now \eqref{36} can be read as $(y_0,y)=(\wh u,\wh\phi)$ being a formal solution of the convergent equation
\beq\label{38}
f'(y)=y_0\cdot f.
\eeq
By Artin's Approximation Theorem \cite{Art68}, \eqref{38} has convergent solutions $(y_0,y)=(u,\phi)$ arbitrarily close to $(\wh u,\wh\phi)\in(\wh\sO^*,\Aut(\wh\sO))$ in the $\mm$-adic topology.
In particular, we can assume that $(u,\phi)\in(\sO^*,\Aut(X))$.
Then $\phi^*f'=u\cdot f$ yields $\phi(D)=D'$ and hence $D\cong D'$.
\end{Proof}

Corollary~\ref{37} and Proposition~\ref{34} now show that the word ``linear'' in the term ``reductive linear free divisor'' is redundant.

\begin{Theorem}\label{32}
Reductive free divisors are linear.
\end{Theorem}

As a generalization of Theorem~\ref{32} we note the following convergent version of Theorem~\ref{17} for quasihomogeneous singularities.

\begin{Theorem}\label{39}
For a quasihomogeneous singularity $D$, Theorem~\ref{17} holds convergently, that is, with $\wh\dl_D$ and $\wh\mm$ replaced by $\dl_D$ and $\mm$.
\end{Theorem}

\begin{Proof}
For a quasihomogeneous $D$, $\dl_{D,0}$ contains a semisimple element with strictly positive eigenvalues.
By \cite[Thm.~5.4]{GS06}, Remark~\ref{30}.\eqref{30a}, and the conjugacy of all Cartan subalgebras \cite[Ch.~III, \S4, Thm.~2]{Ser01}, such an element occurs in $\hh\oplus\dd$ as constructed in Theorem~\ref{17}.
Then $f$ in Theorem~\ref{17}.\eqref{17f} must be a quasihomogeneous polynomial and we can apply Artin's Approximation Theorem as in \eqref{38} to conclude that $f$ is a defining equation of $D$.
\end{Proof}

Finally, we shall generalize Theorem \ref{39} to algebraic singularities using results of Hauser and M\"uller \cite{Mue86,HM89}.
Denote by
\beq\label{45}
\Lambda_D:=\{\phi\in\Aut\sO\mid\phi^*f\in\sO^*\cdot f\}
\eeq
the group of embedded automorphisms of $D$.
Using the projections
\[
\xymat{\A:=\Aut\sO\ar@{->>}[r]^-{\pi_k}&\Aut(\sO/\mm^{k+2})=:\A_k},\quad k\ge0,
\]
we define the truncations
\[
\Lambda_{D,k}:=\pi_k(\Lambda_D),\quad k\ge0.
\]
By \cite[Lem.~1]{Sch06}, each $\Lambda_{D,k}$ is an algebraic group with Lie algebra $\dl_{D,k}$.
Let $L$ be a Levi subgroup in $\Lambda_{D,0}$.
Its Lie algebra is then of the form $\fs_0\oplus\dd_0\subset\dl_{D,0}$ with $\fs_0$ and $\dd_0$ as in Theorem~\ref{17}.\eqref{17c}.
By Theorem~\ref{17} it lifts to $\fs\oplus\dd\subset\wh\dl_D\cap\gl_n(\CC)$ with respect to some basis $\wh\phi$ of $\wh\mm$.
So if we consider $L\subset\GL_n(\CC)\subset\A_k$ with respect to $\wh\phi$, it follows that $L\subset\Lambda_{D,k}\cap\GL_n(\CC)$ for all $k\ge0$ by \cite[Ch.~V, \S13.1, Thm.]{Hum75}.
Thus, $L\subset\Lambda_{\wh D}$ with the latter defined as in \eqref{45}.
Analogous to the proof of \cite[Thm.~1]{Sch06}, one concludes the following result using \cite[Hilfssatz~2]{Mue86} and \cite[Thm.~2']{HM89}.

\begin{Theorem}\label{46}
For an algebraic singularity $D$, Theorem~\ref{17} holds convergently, that is, with $\wh\dl_D$ and $\wh\mm$ replaced by $\dl_D$ and $\mm$.
Moreover, $\fs\oplus\dd$ is the Lie algebra of a Levi subgroup in $\Lambda_D$.
\end{Theorem}

\section{Singular locus and initial Levi factors}\label{22}

A module $V$ over a semisimple Lie algebra $\fs$ splits into weight spaces
\[
V=\sum_{\alpha\in\hh^*}V_\alpha
\]
with respect to a Cartan algebra $\hh$ of $\fs$.
Any $h\in\hh$ acts on $V_\alpha$ through multiplication by $\alpha(h)$.
By \cite[Ch.~VII, Prop.~3]{Ser01}, the weights are (complexifications of) elements of the dual lattice of that generated by the inverse roots of $\fs$.
As such they form the weight diagram $W$ of $\fs$ on $V$ which is finite and invariant under the Weyl group, see \cite[Ch.~VII, \S 4, Rem.~(1)]{Ser01}.

By the Theorem of Malcev--Harish-Chandra \cite[Ch.~III, \S9]{Jac79}, all Levi factors $\fs$ of a Lie algebra $\fg$ are conjugate by an inner automorphism.
Therefore the weight diagram of $\fs$ on a $\fg$-module $V$ is independent of the choice of the Levi factor $\fs$.

\begin{Definition}\label{11}
For a $\fg$-module $V$ we denote by $W=W(\fg,V)$ the weight diagram of a Levi factor of $\fg$ on $V$.
Setting $d(\{\alpha\})=\dim V_\alpha$ for each $\alpha\in W$ and $d(\emptyset)=-\infty$ defines an additive map $d\colon\sP(W)\to\NN$ where $\sP(A)$ denotes the power set of a set $A$.
For any $k\ge3$, let $\C_k(W)\subset\sP(W)$ be the set of all subsets $C\subset W$ for which 
\beq\label{12}
(\underbrace{C+\cdots+C}_l)\cap W=\emptyset
\eeq
for all $l\ge k-1$.
Finally, we set $M_k(\fg,V)=\max d(\C_k(W))$ for $k\ge3$ and $M_k(\fg,V)=-\infty$ otherwise.
\end{Definition}

\begin{Remark}
If $\sll_2(\CC)\cong\fs\subset\fg$ acts nontrivially on $V$ then $M_k(\fg,V)\ge0$ for all $k\ge3$.
\end{Remark}

Recall that the order $\ord f$ of $f$ is defined as the maximal $k$ for which $f\in\mm^k$ and that $\ord D=\ord f$ for each $f$ defining $D$.

\begin{Theorem}\label{13}
$\dim\Sing D\ge M_{\ord D}(\il_D,\mm/\mm^2)$.
\end{Theorem}

\begin{Remark}\label{23}
\begin{asparaenum}

\item\label{23a} In Example~\ref{7} we have $\fs=\sll_2(\CC)$, $k=\ord D=4$, $W=\{-3,-1,1,3\}$, and $M_k(\il_D,\mm/\mm^2)=1<\dim\Sing D$ by choosing the facet $C=\{3\}$ of the convex hull of $W$ corresponding to the $w$-coordinate.
The reason why the estimate of Theorem~\ref{13} is not sharp here is that \eqref{9} is not a coordinate space but contains the $1$-dimensional $w$-axis. 
This will become clear in the proof of Theorem \ref{13}.

\item\label{23b} For $k\ge3$, a weak lower bound for $M_k(\fg,V)$ is the rank of $\fs$.
To see this, pick $\ell\in\hh$ such that $\ell\le1$ on $W$ and $\ell=1$ defines a facet of the convex hull of $W$ and choose $C$ in Definition~\ref{11} as $C=W\cap\{\ell=1\}$. 
Then $\ell(\underbrace{C+\cdots+C}_{l})=\{l\}$ which implies \eqref{12}.
Moreover, $d(C)\ge|C|\ge\dim C+1=\dim W=\rk\fs$.

\item\label{23c} As $M_k(\fg,V)=-\infty$ for $k<3$ the statement of Theorem~\ref{13} is vacuous if $\ord D<3$.
However if $\ord D=1$, $D$ is smooth and both sides of the equality in Theorem~\ref{13} equal $-\infty$.

In case $\ord D=2$, there is a coordinate system $y=y_1,\dots,y_n$ such that 
\[
f(x)=g(y')+y_{k+1}^2+\cdots+y_n^2,\quad y'=y_1,\dots,y_k,
\]
and $D'=\{g=0\}\subset(\CC^k,0)=X'$ is a hypersurface of order $\ord D'\ge3$ with local ring $\sO'=\sO_{X'}$.
So Theorem~\ref{13} can be applied to $D'$ instead of $D$ to get a lower bound for $\Sing D=\Sing D'\subset X'$.
However it is not clear how Levi factors of $\il_D$ relate to those of $\il_{D'}$.

\end{asparaenum}
\end{Remark}

By Remark~\ref{23}.\ref{23b} and \ref{23}.\ref{23c}, Proposition~\ref{6} is a consequence of Theorem~\ref{13}.
In the following we deduce it directly from Theorem~\ref{17}.

\begin{Proof}[Proof of Proposition~\ref{6}]
For a nonquasihomogeneous $D$, $\il_D$ is nilpotent by \cite[Prop.~4]{Sch07}.
We may hence assume that $\eps(f)=f$ for some Euler vector field $\eps\in\dl_D$.
Then any other $\delta\in\dl_D$ may be replaced by $\delta-\frac{\delta(f)}{f}\cdot\eps$ such that $\delta(f)=0$ afterward.
This shows that 
\begin{equation}\label{18}
\dl_D=\sO\cdot\eps\oplus\dl_f,\quad\dl_f:=\Der(-\log f)=\{\theta\in\Der_\CC(\sO)\mid\theta(f)=0\}.
\end{equation}
But $\dl_f$ is isomorphic to the syzygy module of the ideal of partials of $f$.
As $f$ has an isolated critical point, $\p_1(f),\dots,\p_n(f)$ is a regular sequence, and hence 
\begin{equation}\label{19}
\dl_f=\ideal{\p_i(f)\cdot\p_j-\p_j(f)\cdot\p_i\mid 1\le i,j\le n}\subset\Delta^1.
\end{equation}
This means that $\dl_f/\mm\cdot\dl_f\subset\ker\ol\lambda_0$ is nilpotent by \eqref{2}, \eqref{40}, and \eqref{41}.
\end{Proof}

Let $\fs$ be a semisimple Lie subalgebra of $\il_D$ and $\hh$ a Cartan subalgebra of $\fs$.
Note that $\hh$ is Abelian consisting of semisimple elements \cite[Ch.~III, Thm.~3]{Ser01}.  
For any $\hh$-homogeneous $p\in\wh\mm$ we denote its $\hh$-weight by
\[
\wt(p):=(h\mapsto h(p)/h)\in\hh^*.
\]
By Theorem~\ref{17}, we can identify $\fs$ with a linear subalgebra $\fs\subset\wh\dl_D$ and $\wh D$ is defined by an $\hh$-homogeneous $f\in\wh\mm$.

\begin{Lemma}\label{20}
$\wh D$ is defined by some $f\in\wh\mm$ with $\wt(f)=0$.
\end{Lemma}

\begin{Proof}
For $k=\ord f$, $0\ne[f]=f+\mm^{k+1}\in\mm^k/\mm^{k+1}$ generates a $1$-dimensional $\il_D$- and hence $\fs$-module.
By the structure of semisimple Lie algebras \cite[Ch.~VI]{Ser01}, $\fs$ is a sum of Lie subalgebras isomorphic to $\sll_2(\CC)$.
The classification of $\sll_2(\CC)$-modules \cite[Ch.~IV]{Ser01} then implies that $\wt(f)=\wt([f])=0$.
\end{Proof}

\begin{Proof}[Proof of Theorem~\ref{13}]
As explained at the end of Section~\ref{21}, it suffices to prove the claim for $\wh D$ defined by $f$ considered as an element of $\wh\mm$.
Let $\fg=\il_D$, $V=\mm/\mm^2$, $W=W(\fg,V)$, and $k=\ord(f)$. 
Pick $C\in\C_k(W)$ such that $M_k(\fg,V)=d(C)$.
The claim follows if we show that $\Sing\wh D$ contains the linear space defined by all variables $x$ with $\wt(x)\not\in C$.
This means that no partial derivative of $f$ contains a monomial $c$ that is a product of only variables $x$ with $\wt(x)\in C$.
Assume, on the contrary, that such a derivative exists.
Then $\wt(c)\in\underbrace{C+\cdots+C}_l$ where $l:=\deg(c)\ge k-1$ denotes the degree of $c$.
Choosing $f$ with $\wt(f)=0$ as in Lemma~\ref{20}, each partial of $f$ is $\hh$-homogeneous with weight in $-W$ which equals $W$ (by invariance of $W$ under the Weyl group).
In particular, $\wt(c)\in W$ by choice of $c$ in contradiction to \eqref{12}.
\end{Proof}

\bibliographystyle{amsalpha}
\bibliography{inll}

\end{document}